%% file: divrig.tex
\documentclass[a4paper, 12pt]{amsart}
\usepackage{amsfonts,amsmath}
\mathchardef\emptyset="001F

\theoremstyle{plain}

\newtheorem{theorem}{Theorem}[section]
\newtheorem{lemma}[theorem]{Lemma}
\newtheorem{proposition}[theorem]{Proposition}

\newtheorem{remark}[theorem]{Remark}

\theoremstyle{definition}
\theoremstyle{remark}
\numberwithin{equation}{section}

\newcommand{\Om}{\Omega}

\newcommand{\weakst}{\stackrel{\ast}{\rightharpoonup}}

\newcommand{\rank}{\mathrm{rank}}
\renewcommand{\div}{\mathrm{div}}
\newcommand{\Div}{\mathrm{Div}}

\newcommand{\M}{{\mathbb M}}

\newcommand{\hs}{{\mathcal H}}

\newcommand{\R}{{\mathbb R}}

\newcommand{\N}{{\mathbb N}}

\newcommand{\dist}{{\rm dist}}

\newcommand{\dpr}{{\mathcal D}'}

\newcommand{\qu}{{Q}}
\renewcommand{\l}{\mathbb L}

\title[The three divergence free matrix fields problem]{
The three divergence free matrix fields problem}
\author[Mariapia Palombaro  ]{Mariapia Palombaro}
\address[Mariapia Palombaro]{Dipartimento di Matematica, 
Universit\`a ``La Sapienza'', P.le Aldo Moro 2, 00185 Roma, Italy}
\email{palombar@mat.uniroma1.it.}
\author[Marcello Ponsiglione]{Marcello Ponsiglione}
\address[Marcello Ponsiglione]{S.I.S.S.A., Via Beirut 2-4, 34014,
Trieste, Italy}
\email{ponsigli@sissa.it}

\usepackage{graphicx}
\input psfig.sty
\begin{document}

\baselineskip3.15ex
\vskip .3truecm

\begin{abstract}
\small{
We prove that for any connected open set $\Om\subset \R^n$ 
and for any set of matrices 
$K=\{A_1,A_2,A_3\}\subset\M^{m\times n}$, with $m\ge n$ and 
rank$(A_i-A_j)=n$ for $i\neq j$, there is no non-constant solution 
$B\in L^{\infty}(\Om,\M^{m\times n})$, called exact solution, 
to the problem 
\begin{equation*}\label{abs}
\Div B=0 \quad \text{ in }\dpr(\Om,\R^m) \quad 
\text{and} \quad B(x)\in K \text{ a.e. in } \Om\,.
\end{equation*} 
In contrast, A. Garroni and V. Nesi \cite{GN} 
exhibited an example of set $K$ for which the 
above problem admits the so-called approximate solutions. 
We give further examples of this type.  

We also prove non-existence of exact solutions 
when $K$ is an arbitrary set of matrices satisfying 
a certain algebraic condition which is weaker than 
simultaneous diagonalizability.    
 
}
\end{abstract}
\maketitle
\tableofcontents

\section{Introduction}
The problem of characterizing solenoidal matrix 
fields which take values in a finite set of matrices, 
has been recently considered by A. Garroni and V. Nesi. 
This kind of problem is analogous to that on curl free matrix fields 
in which one asks whether a Lipschitz mapping 
using a finite number of gradients exists. 
Here the differential constraint of being the gradient 
of a mapping, and hence a curl free matrix field, 
is replaced by that of being a divergence free matrix field 
(i.e.\hspace{-4pt} a matrix valued function whose rows are divergence free 
in the distributional sense). To describe the problem we begin 
with some definitions.

\noindent
{\it Definition 1. } 
Given two integers $m,n\ge 2$, a set 
of real $m\times n$ matrices $K\subset \M^{m\times n}$ 
and a bounded open set $\Om$ in 
$\R^n$, we say that any $B\in L^{\infty}(\Om,\M^{m\times n})$ 
satisfying 

\begin{equation}\label{pb1}
\begin{cases}
\Div B=0     &  \text{in }\dpr(\Om,\R^m)\,,\\
B(x)\in K    &   \text{a.e. in } \Om\,,\\
B \text{ is non-costant\,,} 
\end{cases}
\end{equation}
\vspace{1 mm}

\noindent
is an exact solution of \eqref{pb1}. 
We say that $K$ is rigid for exact 
solutions if there is no solution to \eqref{pb1}. 
\vspace{1 mm}

\noindent
{\it Definition 2. } 
We say that Problem (\ref{pb1}) admits an approximate solution 
if there exists a uniformly bounded sequence 
$\{B_h\}\subset L^{\infty}(\Om,\M^{m\times n})$ such that

\begin{equation}\label{pb2}
\begin{cases}
\Div B_h\to 0 &\text{in }W^{-1,\infty}(\Om,\R^m), \\
\dist(B_h,K)\to 0 &\text{in  measure},\\
\dist(B_h,A)\not\to 0 &\text{in measure, for every }A\in K.     
\end{cases}
\end{equation}
\vspace{1 mm}

\noindent  
We say that $K$ is rigid for approximate solutions of 
\eqref{pb1} if there is no solution to \eqref{pb2}. 
\vspace{1 mm}

\noindent 
We remark that if $K$ is rigid for approximate solutions and there exists a 
sequence $\{B_h\}$ satisfying the 
first two conditions of (\ref{pb2}), then 
any accumulation point of the sequence $\{B_h\}$ has to be a 
constant matrix in $K$. 
      
Let us briefly describe the situation 
in the context of the ``gradient problem'', that is: 
find $f\in W^{1,\infty}(\Om,\R^m)$ such that $Df\in K$ a.e. in $\Om$ 
and $f$ is not affine.  
It is well-known in this setting, that 
a sufficient condition for a set $K=\{A_1,\dots,A_N\}\in\M^{m\times n}$, 
$N\le 3$, to be rigid is that rank$(A_i-A_j)>1$ for $i\neq j$. 
The condition  rank$(A_i-A_j)=1$ is called rank-1 connectedness. 
J.~ M. Ball and R. D. James studied in detail the case $N=2$ and 
proved, under the latter assumption, a rigidity result both for 
exact and for approximate solutions (see \cite{Ball87}).  
For $N=3$ rigidity still holds. The following theorem is due to 
\v Sver\'ak and will be used later.
\begin{theorem}\label{teosv}{\rm (V. \v Sver\'ak, \cite{Sver92})}. 
Let $\Om\subset\R^n$ be an open connected set and 
let $K=\{A_1,\,A_2,\,A_3\}\subset \M^{m\times n}$, with 
{\rm rank}$(A_i-A_j)> 1$ for $i\neq j$. 
If $f\in W^{1,\infty}(\Om,\R^m)$ satisfies $Df\in K$ a.e., then 
$Df$ is constant.

\noindent
Let $p>2$ and let $f_h\rightharpoonup f$ in $W^{1,p}(\Om,\R^m)$. 
If $\dist(Df_h,K)\rightarrow 0$ in measure, then $Df_h\rightarrow A_i$ 
in measure, for some $i\in\{1,2,3\}$. 
\end{theorem}

\noindent
The previous result, specialized to the case $m=n=2$, 
will be a crucial tool in the 
proof of Theorem \ref{main}.

\noindent
For completeness let us also recall that 
for $N=4$ rigidity still holds for exact solutions 
(see \cite{Ch02}) but it can fail for approximate ones 
and a suitable choice of $\{A_1,A_2,A_3,A_4\}$ (see \cite{Tar87}, \cite{Tar93}). 
The case $N=5$ is nicely illustrated in \cite{Kir} by a non-rigid five 
point configuration without any rank-1 connection.

\begin{remark}\label{equi}{\rm All the previous works provide results also for the 
``divergence problem'' when the working space is $\M^{m\times 2}$, 
since any set of solenoidal matrix fields defines a set of gradients 
via right-multiplication by 
$J:=\left(
\begin{array}{rr}
0 & -1 \\
1 & 0
\end{array}\right)$.}
\end{remark} 
The problems are no longer equivalent if $m\geq n>2$. 
The right notion of connectedness which comes into play, in this case, 
is that of the so called rank-$(n-1)$ connectedness. Given  
$A_1,A_2\in\M^{m\times n}$, with $m\ge n$ and $\rank(A_1-A_2)\le n-1$, 
one can construct solenoidal matrix fields which take 
both the values $A_1$ and $A_2$ on a set of positive measure 
(indeed one can check that simple laminates work). In contrast, 
if $\rank(A_1-A_2)=n$, the following rigidity result holds.    
\begin{proposition}\label{teogn1}{\rm (A. Garroni, V. Nesi, \cite{GN})}. 
Let $\Om$ be an open and connected set in $\R^n$. 
Let $A_1,A_2\in\M^{m\times n}$, with $m\ge n$ and $\rank(A_1-A_2)=n$. 
Let $B:\Om\rightarrow \{A_1,A_2\}$ be a measurable function satisfying 
$\Div B=0$ in $\dpr(\Om,\R^m)$. Then B is constant.
\end{proposition}

\noindent
The case when $K$ is made of two matrices has been completely solved 
and given a negative answer also for what concerns approximate 
solutions. 
The following proposition establishes rigidity for approximate 
solutions under the hypothesis of
rank-$(n-1)$ disconnectedness.
\begin{proposition}\label{teogn2}{\rm (A. Garroni, V.  Nesi, \cite{GN})}.
Let $\Om$ be a bounded open and connected set in $\R^n$, 
and let $K=\{A_1,A_2\}\subset\M^{m\times n}$, $m\geq n\ge 1$, be 
such that $\rank(A_1-A_2)=n$. Let $B_h$ be a sequence weakly convergent 
to $B$ in $L^p(\Om,M^{m\times n})$, with $p>1$, such that
$$\Div B_h\rightarrow 0\quad \text{ strongly in } W^{-1,p}(\Om,\R^m)$$
and
$$\dist(B_h,K)\rightarrow 0\quad \text{ in measure}\,.$$
Then
$$B_h\rightarrow A_1 \quad \text{ or }\quad B_h\rightarrow A_2 
\quad \text{ in measure}\,.$$
\end{proposition}

\noindent
So far everything seems to be parallel the ``gradient problem'', 
but the case when $K$ consists of three matrices turns out to be 
different.  
Indeed one can construct a sequence of matrix fields which 
are divergence free and whose distance from the set $K$ approaches 
zero. In other words, approximate solutions exist for a suitable 
choice of $\{A_1,\, A_2,\, A_3\}$. The following result 
clarifies the situation.
\begin{lemma}\label{app}{\rm (A. Garroni, V. Nesi, \cite{GN})}. 
Given $m\geq n\geq 3$, there exist three pairwise rank-n connected 
$m\times n$ matrices $A_1$, $A_2$, $A_3$, and there exists a 
sequence $B_h\in L^{\infty}(\Omega,\M^{m\times n})$ such that 
setting $K=\{A_1,\, A_2,\, A_3\}$, one has
\begin{align}\label{goal1}
&{\rm dist}(B_h,K)\to 0 \quad\quad \text{ strongly in }L^p(\Om),\,
\forall p \ge 1,\\
&\label{goal2} \Div B_h \to 0 \quad\quad \text{ strongly in }W^{-1,p}(\Om;\R^m),\,
 \forall p \ge 1,
\end{align}
and $B_h \weakst B$ in $L^\infty$, with $B\neq A_i$ for any $i=1,\,2,\,3$.
\end{lemma}

\begin{remark}{\rm (A. Garroni, V. Nesi, \cite{GN}).  
In Lemma \ref{app} one can achieve the stronger requirement $\Div B_h=0$ 
rather then \eqref{goal2}, by suitably projecting the fields $B_h$ onto 
Divergence-free matrix fields.} 
\end{remark}

\noindent
The explicit formula for $A_1$, $A_2$, $A_3$, can be found in the 
last section of this paper (see Remark \ref{matgn}). 

Next, we state the main theorem of our paper, namely a rigidity 
result for three-valued matrix fields under the assumption of 
rank-$(n-1)$ disconnectedness. 
\begin{theorem}\label{main}
Let $\Om\subset\R^n$ be a connected open set and    
let $K=\{A_1$, $A_2$, $A_3\}\subset \M^{m\times n}$,
with  $m\ge n$ and $\rank(A_i - A_j)=n$ for $i\neq j$. 
If $B:\Om\to K$ is a measurable function 
satisfying $\Div B=0$ in $\dpr(\Om,\R^m)$, then $B$ is constant. 
\end{theorem}
\vspace{1em}
In this paper we will mainly deal with the problem of non-existence of 
exact solutions. In addition to the previous theorem, we prove a   
rigidity result for a particular class of matrix fields taking an 
arbitrary number of values.  
This is the precise statement. 
\begin{theorem}\label{d}
Let $\Om \subset \R^n$ be a connected open set and let 
$K\subset\M^{m\times n}$ be bounded, with $m\ge n$ and rank$(A_i-A_j)=n$ for 
every $A_i,A_j\in K$, with $i\neq j$. 
Suppose that $K$ satisfies the following condition:

\noindent
there exist $n-1$ independent hyperplanes 
$\pi_1,\dots,\pi_{n-1}$ in $\R^n$, and $n-1$ vector subspaces 
$\tau_1,\dots,\tau_{n-1}$ of dimension $n-1$ in $\R^m$, such that  
\begin{equation}\label{conditio}
A_i: \pi_r\to A_i(\pi_r)\subseteq \tau_r\,, \text{ for every } A_i\in K 
\text{ and } r=1,\dots n-1\,. 
\end{equation}
Then every measurable matrix field 
$B:\Om \to K$ satisfying 
$\Div B=0$ in $\dpr(\Om,\R^m)$ is constant.
\end{theorem}

\noindent
We remark that Theorem \ref{d}, specialized to the case $m=n=2$, reduces 
itself to a result which is well-known in the setting of 
the gradient problem ( see \cite{Ch02}, Lemma 5). 
An easy corollary of this theorem is that a set of simultaneously diagonalizable 
rank-$n$ connected matrices is rigid for exact solutions. 

The study of more general linear differential constraints on 
the matrix field $B$ is just beginning. The interested reader is 
referred to \cite{marco} for results in this direction. 

\vspace{1em}
The plan of the present paper is as follows. 

\vspace{1em}
In Section \ref{prel}, we present an algebraic argument which implies 
that the right dimension 
to study the problem is $n\times n$ (see Lemma \ref{quad}). 
Next we remark that the condition of being divergence-free 
is invariant under any orthogonal change of variables (Remark \ref{rem}). 
Using this invariance, in order to prove Theorem \ref{main} 
it is enough to consider a very special situation. 
This kind of argument does not work for an arbitrary number of matrices
and one does in fact expect that rigidity fails for a sufficiently large 
number of them. 
Yet, under the assumptions of Theorem \ref{d}, one can prove that 
rigidity still holds for  
an arbitrary number of values and actually even for a continuum 
of them. For the reader's convenience we give, 
in Lemma \ref{gg}, the Gauss-Green formula for $L^{\infty}$ fields, 
which will be the main ingredient in the proof of this result.    

\vspace{1em}
In Section \ref{mr}, we give the proofs of Theorem \ref{main} 
and Theorem \ref{d}.   

\vspace{1em}
The final section departs from the main focus of the paper. 
Indeed, in the spirit of Lemma \ref{app}, we address the problem of finding 
approximate solutions to the ``three divergence problem''.  
More precisely, we show 
that the construction used by Garroni and Nesi actually applies to a 
larger class of sets $K$. Theorem \ref{fine} gives a characterization 
of all such $K$'s, which turn out to be non-rigid for approximate 
solutions.

\begin{theorem}\label{fine}
For every $q_1$, $q_2$, $q_3\in(0,1)$, let $A\in\M^{3\times 3}$  
be defined as follows 
\begin{equation*}A=\frac{1}{q_3}\left[\left(1-\prod\limits_{i=1}^3(1-q_i)
\right)
G^{-1}
\left(\begin{array}{ccc}
\lambda_1 & 0         & 0         \\
0         & \lambda_2 & 0         \\
0         & 0         & \lambda_3 \\
\end{array}\right)
G-
q_2(1-q_3)I
\right]\,,
\end{equation*}
where ${\displaystyle 
\lambda_1=0\,,\lambda_2=1/(1-q_1)\,, 
\lambda_3=q_2/(q_1+q_2-q_1q_2)}$ , and 
$G$ is an arbitrary matrix in $GL(3)$. 
Then, for every ${\displaystyle M\in\M^{3\times 3}}$ and $N\in GL(3)$, 
the set 
$$K=\{ M\,,N+M\,,NA+M \}$$
is non-rigid for approximate solutions.
\end{theorem}

\section{Preliminaries}\label{prel}
\noindent
In this Section we set some notations and a few preliminary results 
needed in the proof of the main results of Section 3.
\\ Throughout this paper  $\Om$ is  an open connected subset of $\R^n$. 
We denote
by $\M^{m\times n}$ the set of the real $m\times n$ matrices; $0$ and $I$
will indicate the zero matrix and the identity matrix in 
$\M^{n\times n}$ respectively. Left-multiplication of $A$ 
times a vector $v\in\R^n$ is denoted  by $A\cdot v $.  
The symbol $\langle v ,w\rangle$ denotes the standard inner product in $\R^n$.
   
\par
For every measurable subset $E$ of $\R^n$, $|E|$ is the
$n$-dimensional Lebesgue measure of $E$ while, for $s\in\R^{+}$, 
$\hs^s(E)$ is its $s$-dimensional Hausdorff measure. 
\\ Given a function $f\in L^1(\Om,\R)$, we say that $x\in\Om$ is a Lebesgue 
point for $f$, and that
$\lambda(x) \in \R$ is the Lebesgue value of $f$ at $x$,
if 
$$
\lim_{r\to 0}-\hspace{-13 pt}\int_{B_r(x)} |f(y)-\lambda(x)| \, dy =0,
$$
where $B_r(x)$ is the open ball of radius $r$ and center $x$ and 
the symbol $-\hspace{-10 pt}\int_{B_r(x)}$ stands for 
$\frac{1}{|B_r(x)|}\int_{B_r(x)}$. 
This definition extends in the obvious way to vector valued functions.\\  
It is well-known that the set of Lebesgue points for $f$, which 
from now on we will denote by $\l(f)$, has full measure in $\Om$. 
For every $k\in\N$ and $f\in L^1(\Om,\R^k)$, 
we will denote by $\tilde{f}$ a Lebesgue
representative of $f$ (i.e. $\tilde{f}(x)=\lambda(x)$ for 
every $x\in\l(f)$), 
so that $\tilde f$  coincides with  $f$ a.e. in $\Om$.  
For more details we refer the reader to \cite{EG}. 

Recall that a vector field $f\in L^\infty(\Om,\R^n)$ is 
said to be divergence free if for every $\varphi\in C_0^{\infty}(\Om)$
$$
{\displaystyle \int_{\Om}\langle f(x),\nabla\varphi(x)\rangle\,dx=0}.
$$ 
For the reader's convenience, we prove a Gauss-Green formula in the
particular setting of our problem (much more general results can be 
found in \cite{anz}). 
In the sequel the symbol $\nu(x)$ will denote the outward normal to a 
given surface at the point $x$. 

\begin{lemma}\label{gg}
Let $f\in L^\infty(\Om,\R^n)$ be a divergence free vector field, and let 
$U\subset \overline{U} \subset \Om$ be open with Lipschitz boundary. 
Suppose that 
$\hs^{n-1}(\l (f) \cap \partial U)=\hs^{n-1}(\partial U)$. 
Then 
$$
\int_{\partial U}\langle \tilde f (s), \nu(s)\rangle\, d\hs^{n-1}(s) = 0\,.
$$
\end{lemma}

\begin{proof}
Consider a sequence $\{\rho_n\}$ of mollifiers and set $f_n:=f*\rho_n$. 
We have $$\div f_n=\div f*\rho_n=0.$$ The standard 
Gauss-Green formula for smooth functions on a Lipschitz domain yields
\begin{equation}\label{for1}
\int_{\partial U} \langle f_n(s),\nu(s)\rangle\, d\hs^{n-1}(s) = 0.
\end{equation}
It is easy to see that  
$f_n(x)\to \tilde f(x)$ for all $x \in \partial U\cap\l(f) $. 
Passing to the limit in \eqref{for1} and using 
the dominated convergence Theorem, we get
$$
\int_{\partial U} \langle\tilde f(s),\nu(s)\rangle\, d\hs^{n-1}(s) = 0.
$$
\end{proof}

\begin{remark}\label{rem}{\rm
Let $B\in L^{\infty}(\Om,\M^{m\times n})$ be a divergence free matrix field and
let $R$ be an orthogonal matrix in $\M^{n\times n}$.
Using a convolution argument as in Lemma \ref{gg} and the 
classical chain rule formula for smooth functions, 
one can check that for every $C\in \M^{m\times n}$ and $F\in\M^{n\times m}$, 
the matrix field $\widehat{B}:
\{y\in \R^n\,|\, y=R^Tx,\,x\in\Om\}   \to \M^{n\times n}$ defined by
\begin{equation}\label{base}
\widehat{B}(y):= R^T F \big(B(R y) + C\big) R
\end{equation}
is divergence free.} 
\end{remark}
\noindent
The next lemma shows that, given a set $K\subset\M^{m\times n}$, the property
of rank-$(n-1)$ disconnectedness is preserved under left multiplication by 
suitable matrices in $\M^{n\times m}$. 

\begin{lemma}\label{cor1}
Let $K\subset \M^{m\times n}$ be at most countable,  
with $m>n$ and $\rank(A_i - A_j)=n$ for every $A_i$, $A_j\in K$, with $i\neq j$. 
Then there exists $F\in\M^{n\times m}$ such that 
$\rank(FA_i-FA_j)=n$ for every $A_i$, $A_j\in K$ with $i\neq j$.
\end{lemma}

\begin{proof}
Since the image of $A_i - A_j$ is a $n$-dimensional subspace of $\R^m$,  
we can always find a linear operator $F:\R^m\to\R^n$ such that 
${\rm Ker}(F(A_i-A_j))={\rm Ker}(F) \cap {\rm Im}(A_i-A_j) =\{0\}$ for every
$i\neq j$. Then $\rank(FA_i-FA_j)=n$ for every $A_i$, $A_j\in K$ 
with $i\neq j$. 
\end{proof}
\noindent
The previous results show that, as long as we consider discrete 
valued matrix fields, we can always set the problem in the 
space of square matrices. This is the claim of the next Lemma. 

\begin{lemma}\label{quad}
Let $K\subset \M^{m\times n}$ be at most countable, with $m>n$           
and $\rank(A_i-A_j)=n$ for every $A_i,A_j\in K$, with $i\neq j$. 
Suppose we are given a divergence free matrix field 
$B:\Om\to K$. 
Then there exist a divergence free matrix field 
$B':\Om\to K'$, where $K'\subset\M^{n\times n}$,  
${\rm card}(K')={\rm card}(K)$,  
and $\rank(A_i'-A_j')=n$ for every $A_i',A_j'\in K'$, with $i\neq j$.
\end{lemma}
\begin{proof}
Apply Lemma \ref{cor1} 
and use Remark \ref{rem} with $R=I$ and $C=0$.
\end{proof}

\section{Proofs of the main results}\label{mr}
\noindent
In this section we give the proofs of Theorem \ref{main} and 
Theorem \ref{d}. 

\begin{proof}[{\bf Proof of Theorem \ref{main}}]
By Lemma \eqref{quad} it is enough to consider  
the case $m=n>2$. 
Moreover, due to the local character of our problem, we can make any 
convenient change of variables. Hence, 
as it is customary in this kind of problems, we begin with some
reductions to special cases. 
We use Remark \ref{rem} choosing $F=(A_2 - A_1)^{-1}$
and $C=-A_1$. In this way we can assume that $A_1=0$ and $A_2=I$.
Moreover, since for any linear operator in $\R^n$ there exists a 
two-dimensional invariant subspace, we can choose the orthogonal 
matrix $R$ in \eqref{base} so that  $A_3$ is of  the form
\begin{equation}\label{partf}
A_3 := A =\left(
\begin{array}{ccccc}
a_{11} & a_{12} & 0 & \cdots & 0\\
a_{21} & a_{22} & 0 & \cdots & 0 \\
a_{31}  & a_{32} & a_{33} & \cdots   &  a_{3n} \\
\vdots  & \vdots & \vdots & \ddots & \vdots \\
a_{n1}  & a_{n2}& a_{n3} & \cdots & a_{nn}
\end{array}
\right).
\end{equation}
and rank$(A)=$ rank$(A-I)=n$.
By assumptions we can write 
$$B=\chi_{E_1}0+\chi_{E_2}I +\chi_{E_3}A\,,$$ 
where $E_i$ are disjoint measurable sets and 
$E_1\cup E_2\cup E_3 =\Omega$.  
Next remark that, due to \eqref{partf}, the  
first two equations of $\Div B=0$ only 
involve derivatives with respect to directions $e_1$ and $e_2$. 
Roughly speaking, the idea is to use these informations   
to conclude that $B$ 
does not depend on the variables $x_1$ and $x_2$.
This allows us to make a ``projection'' of the original problem 
into a lower dimensional space (actually a section of $\Om$) and to 
proceed by induction on the dimension~$n$, 
using, at the final step, that by Theorem \ref{teosv} and Remark \ref{equi}, 
rigidity for exact solutions holds for $n=2$. 

Let us proceed with the formal proof. 
Set $n>3$ and suppose that rigidity holds in $\M^{(n-1)\times (n-1)}$. 
We want to prove that $B$ is constant. 
For every $x\in \Om$, let $Q$ be a open coordinate cube 
centered in $x$ and such that 
$\overline{Q}\subset \Om$.
Without loss of generality we assume that $x=0$,
so that $\qu = (-l,l)^n$ for some positive $l$.

Let $\{\rho_i\}$ be a sequence of mollifiers and set $B_i:= \rho_i * B$.
For $i$ large enough we have that $B_i$ is well defined 
in $Q$, $B_i \in C^{\infty}(Q,\M^{n\times n})$ and 
 $\Div B_i = 0$ in $Q$ in the classical sense.
Then fix $\overline{x}:=(\overline{x}_3, \cdots , \overline{x}_n) \in (-l,l)^{n-2}$ 
such that for $\hs^2$-a.e.
$(x_1,x_2) \in (-l,l)^{2}$, $(x_1,x_2,\overline{x})\in \l(B)$. 
By Fubini's Theorem, this is possible for $\hs^{n-2}$-a.e. 
$\overline{x}\in (-l,l)^{n-2}$, 
since, as already remarked, $\l(B)$ has full 
measure in $\Om$.  
Now consider the field $\underline{B}: (-l,l)^2 \to \M^{2\times 2}$ defined by

\begin{equation}
\underline{B}(x_1,x_2):= 
\left(
\begin{array}{cc}
\tilde{b}_{11} & \tilde{b}_{12}\\
\tilde{b}_{21} & \tilde{b}_{22}
\end{array}
\right) (x_1,x_2,\overline{x}),
\end{equation} 
where $\tilde{b}_{kl}$ is a Lebesgue representative of the 
$kl$-entry of the matrix $B$. 
Then let $b^i_{kl}$ be the $kl$-entry of $B_i$ and set

\begin{equation}
\underline{B_i}(x_1,x_2):= 
\left(
\begin{array}{cc}
b^i_{11} & b^i_{12}\\
b^i_{21} & b^i_{22}
\end{array}
\right) (x_1,x_2,\overline{x}).
\end{equation} 

\noindent 
Since $\underline{B}_i \in C^{\infty}((-l,l)^2,\M^{2\times 2})$ and 
 Div $\underline{B}_i = 0$ in $(-l,l)^2$, we have
\begin{equation}\label{fd}
\int_{(-l,l)^2} \underline{B}_i \cdot\nabla \varphi \, dx_1 dx_2 =0
\quad\quad \forall \varphi
\in C^\infty_0(-l,l)^2.
\end{equation}

\noindent
Moreover $\underline{B}_i(x_1,x_2)$ converges to $\underline{B}(x_1,x_2)$ at every
$(x_1,x_2)\in (-l,l)^2$ such that $(x_1,x_2,\overline{x})$ $\in \l(B)$. 
Then passing to the limit for 
$i\to \infty$ in \eqref{fd} and using the dominated convergence 
Theorem we get
 
\begin{equation}\label{fd2}
\int_{(-l,l)^2} \underline{B}(x_1,x_2)\cdot \nabla \varphi(x_1,x_2) 
\, dx_1 dx_2 =0\quad\quad \forall \varphi
\in C^\infty_0(-l,l)^2.
\end{equation}

\noindent
By Theorem \ref{teosv} we have that $\underline{B}$ is constant and hence 
$\widetilde{B}$ is constant on the section 
$(x_1,x_2,\overline{x})$. 
Since this is true for   $\hs^{n-2}$-a.e. $\overline{x}\in (-l,l)^{n-2}$, 
we deduce that 
$\widetilde{B}$ does not depend on $(x_1,x_2)$ in $Q$.
In particular we have 

\begin{equation}\label{derzer}
\frac{\partial \chi_{E_i}}{\partial x_1} = 0 \quad\quad \text{ in }\mathcal{D}'(\qu).
\end{equation}
From \eqref{derzer}, it is easy to see that there exist three measurable sets
 $E_1'$,  $E_2'$,  $E_3'$ in $(-l,l)^{n-1}$ such that 
\begin{equation}\label{carte}
E_i= (-l,l) \times E_i'  \quad \text{ a.e.}
\end{equation}
Now call $0'$, $I'$, $A'$ the $n\times (n-1)$-minors of the 
matrices $0$, $I$, $A$ respectively, obtained 
by eliminating the first column of each matrix. 
Notice that rank$(A')=$ rank$(A'-I')=n-1$. 
Then set  
$$
B':= 0' \chi_{E_1'}+ I'\chi_{E_2'} + A' \chi_{E_3'}\,.
$$
Let us emphasize that $B'$ leaves in a space of dimension $n-1$. 
Combining  \eqref{derzer} with the equation 
$\Div B=0$ in  $\mathcal D'(Q,\R^n)$, one concludes that $B'$ satisfies 
\begin{equation}\label{rettan}
\Div B' = 0 \quad  \text{ in }  \mathcal D'((-l,l)^{n-1},\R^n)\,.
\end{equation}  
By Lemma \ref{cor1}, there exists  $F\in M^{(n-1)\times n}$ 
such that 
\begin{equation}\label{nuovopb}
\rank(FI')=\rank(FA'-FI')=n-1\,. 
\end{equation}
Then set 
\begin{align*}
& A_2':=FI'\,,\\
& A_3':=FA'\,,\\ 
& B_{(n-1)}:= 0 \chi_{E_1'}+ A_2'\chi_{E_2'} + A_3' \chi_{E_3'}\,.
\end{align*}   
By  \eqref{rettan} and \eqref{nuovopb}, it follows that
$B_{(n-1)}$ is an exact solution of the problem 
\begin{equation}\label{nuovaeq}
\Div B_{(n-1)} = 0 \quad \text{ in } \mathcal D'((-l,l)^{n-1},\R^{n-1}) \,,  
\quad B_{(n-1)}\in\{0,A_2',A_3'\}\subset\M^{(n-1)\times (n-1)}\,.
\end{equation}
\noindent
By the inductive assumption, problem \eqref{nuovaeq} is rigid. 
Then $|(-l,l)^{n-1}| = |E_i'|$ 
for some $i$ and hence, 
by \eqref{carte}, $|Q| = |E_i|$. By the arbitrariness of $Q$ 
we conclude that $B$ is constant.
\end{proof}

\begin{remark}\label{remcil}{\rm
Before giving the proof of Theorem \ref{d}, we make some 
considerations about its assumption \eqref{conditio}. 
We want to show that, even in this case, we can reduce to a very special 
situation. Indeed, condition \eqref{conditio} implies that 
${\displaystyle {\rm Im}(A_i)={\rm Im}(A_j)}$ for every $A_i,A_j\in K$.
We can thus apply the same argument as in 
Lemma \ref{quad}, although no requirement is made 
on the cardinality of the set $K$. More precisely, we fix any two of 
the matrices in $K$, say $A_1$ and $A_2$, and choose a matrix 
$F\in\M^{n\times m}$ such that $F(A_1-A_2)=I$. For such an $F$, we have 
that $FA_i\in \M^{n\times n}$ and rank$(FA_i- FA_j)=n$ for every $A_i,A_j\in K$. 
Moreover, the hyperplanes $\pi_1,\dots,\pi_{n-1}$ are preserved under 
the action of every $FA_i$. 
}
\end{remark}
We are now ready to prove Theorem \ref{d}.

\begin{proof}[{\bf Proof of Theorem \ref{d}}]
By Remark \ref{rem} and Remark \ref{remcil}, we can assume 
that $K\subset \M^{n\times n}$ and $A_i(\pi_r)\subseteq\pi_r$ 
for every $A_i\in K$ and $r=1,\dots,n-1$. 
For every $r$, let $v_r$ be the unit vector orthogonal 
to $\pi_r$. 
We want to prove that $B$ does not depend on any of the 
directions $v_r$, which are independent by assumptions. 
Then $B$ would only depend on one direction, but the 
condition of being divergence free will imply that $B$ is constant.  
We will only check the statement 
for one vector $v_r$. 
Choose any of the vectors $v_r$ 
and a real number $q$, and by contradiction assume that 
there exist two points $P$, $Q\in\l(B)$, with $P-Q = q v_r$ 
and such that $\tilde{B}(P)\neq \tilde{B}(Q)$; 
for instance suppose that  $\tilde{B}(P)=A_1$ 
and $\tilde{B}(Q)=A_2$.

Let us briefly digress to explain the idea of the proof 
in an informal way. 
We want to apply the Gauss-Green formula in a small 
cylinder with axis parallel to $v_r$ and bases centered in
$P$ and $Q$ respectively. Computing the flux of $B$ through 
the boundary of the cylinder we will check that its contribution 
through the bases can never compensate the 
contribution through the lateral boundary, the former 
being a vector parallel to $(A_1-A_2)v_r$, the latter 
belonging to the hyperplane $\pi_r$. 
This idea, however, will require some technical efforts, 
as $B$ may be non-constant on the bases. What we do, in fact, 
is to consider a sequence of cylinders with vanishing radii. 
The contradiction will arise for a sufficiently small radius.

Now we continue the formal proof. 
To simplify the notations we will assume that
$P=(0,\dots,0)$ and $v_r = e_n$, so that $Q=(0,\dots ,0,q)$.  

\noindent     
Let  $\rho (x):=(x_1^2+\dots+x_{n-1}^2)^{1/2}$ and, for every 
$r\in \R^+$, set \\[-5pt]
$$
C(r):= \{x \in \R^n : \quad 0\le \rho(x) \le r, \, 
0\le x_n \le q \}.
$$\\[-5pt] 
Since $\widetilde{B}(P)=A_1$ and $\widetilde{B}(Q)=A_2$, 
we can find a cylinder on whose bases the mean value 
of $B$ is $A_1'$ and $A_2'$ respectively, 
with $|A_1-A_1'|$ and $|A_2-A_2'|$ arbitrarily small 
(this can be checked by using Fubini's theorem). 
More precisely, for every given  $\delta>0$, 
we can find a radius $r_\delta \in \R^+$ and a vector 
$ w_\delta \in \R^n$ such that, setting\\[-10pt]
$$C_\delta := C(r_\delta) + w_\delta$$
and denoting by $D^1_\delta$ and $ D^2_\delta$ the bases of $C_\delta$ and by 
$L_\delta$ its lateral boundary (see Figure 1),  
the following hold:\\[-10pt]
\begin{align}\label{teofub}
\nonumber & C_\delta \subset \Om\,,\\
\nonumber & \hs^{n-1}(\partial C_{\delta}\cap\l(B))=\hs^{n-1}(\partial C_\delta)\,,\\
          &-\hspace{-13 pt}\int_{D^1_\delta}|\widetilde{B}(s)-A_1|d\hs^{n-1}(s) + 
           -\hspace{-13 pt}\int_{D^2_\delta}|\widetilde{B}(s)-A_2|d\hs^{n-1}(s)< 
           \delta\,.  
\end{align}

\begin{center}\label{figura}
\begin{figure}\input{div7.pstex_t}\caption{The cylinder $C_\delta$}
\end{figure}
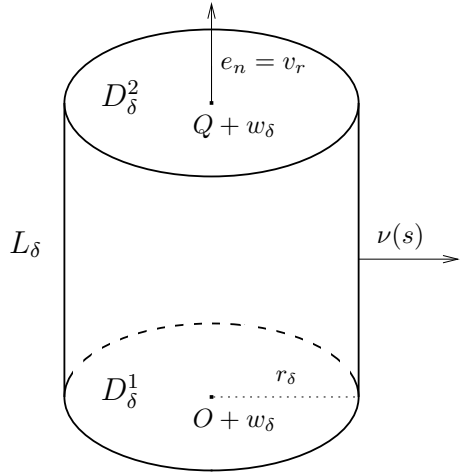
\end{center}

\noindent
By Lemma \ref{gg} we get

\begin{equation}\label{dize}
\int_{\partial  C_\delta} \widetilde{B}(s)\cdot \nu(s) \, d\hs^{n-1}(s) = 0\,.
\end{equation}

\noindent 
We write \eqref{dize} as the sum of 
three contributions as follows
   
\begin{align}\label{dize1}
\nonumber
&\int_{\partial  C_\delta} \widetilde{B}(s)\cdot \nu(s) \, d\hs^{n-1}(s) =
\int_{D^1_\delta} \widetilde{B}(s)\cdot (-e_n) \, d\hs^{n-1}(s) +\\
&\int_{D^2_\delta} \widetilde{B}(s)\cdot e_n \, d\hs^{n-1}(s) +
\int_{L_\delta } \widetilde{B}(s)\cdot \nu(s) \, d\hs^{n-1}(s)\,,
\end{align}

\noindent
where $\nu(s)=\frac{1}{\rho(s)}(s_1,\dots,s_{n-1},0)$ on $L_\delta$. 
On the other hand we have
\begin{align}\label{dize2}
\nonumber
\int_{D^1_\delta}\widetilde{B}(s)\cdot (-e_n) \, d\hs^{n-1}(s)& =
\hs^{n-1}(D^1_{\delta})\, A_1^{\delta}\cdot(-e_n)
 \,,\\ 
\int_{D^2_\delta} \widetilde{B}(s)\cdot e_n \, d\hs^{n-1}(s)& =
\hs^{n-1}(D^1_{\delta})\, A_2^{\delta}\cdot e_n\,,
\end{align}
where $A_1^{\delta}$, $A_2^{\delta}\in\M^{n\times n}$ and, 
by \eqref{teofub}, are such that 
$| A_1^{\delta}- A_1| +| A_2^{\delta}- A_2|<\delta $. 
Then we set 
$$u_\delta := \int_{L_\delta } \widetilde{B}(s)\cdot \nu(s) \, d\hs^{n-1}(s)\,.$$
Dividing the right hand side in \eqref{dize1} by $\hs^{n-1}(D^1_{\delta})$ 
and using \eqref{dize} and \eqref{dize2}, we obtain

\begin{equation}\label{fin}
(A_2^{\delta}-A_1^{\delta}) e_n + \frac{u_\delta}{\hs^{n-1}(D^1_{\delta})} = 0 \,.
\end{equation}
Now recall that $e_n$ is orthogonal to $\pi_r$, rank$(A_2-A_1)=n$ and 
$(A_2-A_1)(\pi_r)=\pi_r$. It follows that $(A_2-A_1)\cdot e_n\notin\pi_r$, and 
hence $(A_2^{\delta}-A_1^{\delta})\cdot e_n\notin\pi_r$ 
for $\delta$ small enough.
On the other hand we have that $u_\delta\in\pi_r$, 
the hyperplane $\pi_r$ being preserved under the action of every $A_i\in K$.  
Then, for sufficiently small $\delta$, \eqref{fin} gives a contradiction. 
\end{proof}

\begin{remark}{\rm 
The assumption of boundedness required for the set $K$, in Theorem \ref{d}, 
can be actually removed.  
The previous proof, indeed, can be adapted to $L^1_{{\rm loc}}$ matrix fields 
by suitable modifications.} 
\end{remark}

\begin{remark}
{\rm Theorem \ref{d}, specialized to the case $m=n=2$,  
gives a rigidity result which was already known 
in the setting of the ``gradient problem'' 
(see \cite{Ch02}, Lemma 5).   
}
\end{remark}

\section{Approximate solutions}\label{ult}
\noindent
In this section we give a refinement of Lemma 
\ref{app} in \cite{GN} in which the authors give an explicit example of a set 
which is non-rigid for approximate solutions. 
Their construction is set in $\M^{3\times 3}$, but it 
can be extended to the case $m,n\ge 3$ by slight modifications. 
It actually provides an algorithm  (similar to that of Tartar's for 
the gradients, \cite{Tar93}), which allows us 
to find approximate solutions for a large class of sets $K$. 
Similar constructions can be found in the works of several authors, 
see \cite{au}, \cite{cas}, \cite{nesi}, \cite{sch}, \cite{talb}. 
The particular case here resembles the construction in \cite{nesi}. 
The key point is the following lemma. 
\begin{lemma}\label{propo1}
Let $K=\{A_1,\, A_2,\, A_3\}\subset \M^{3\times 3}$ be a set 
of pairwise rank-3 connected matrices. If there exist three matrices 
$S_1$, $S_2$, $S_3\in \M^{3\times 3}$ which satisfy the conditions 
\begin{align}
\label{cond1}
&\det(A_i-S_i)=0, \quad \text{for }\, i=1,2,3\,, \\
\label{cond2}
&S_i=q_{i-1}A_{i-1}+(1-q_{i-1})S_{i-1} \quad \mod 3\,, \quad \text{for } \,i=1,2,3\,,
\end{align}
for some $q_i \in (0,1)$, then the set $K$ is non-rigid for approximate 
solutions.
\end{lemma}

The proof of this lemma relies on works on multiple scales. 
We refer to \cite{br1} and \cite{br2} for a general treatment 
and to \cite{GN} and \cite{Tar93} for the case of interest here. 
We simply remark that condition (\ref{cond1}) is used  
to laminate  $A_i$ and $S_i$ in some direction belonging to 
Ker$(A_i-S_i)$, while \eqref{cond2} is used to construct 
a sequence which ``approaches'' the set $K$ in 
the sense of Lemma \ref{app}.    
This is the strategy used in \cite{GN} where the 
authors make an explicit choice of the matrices $A_i$ and $S_i$.  
Theorem \ref{fine}, given in the introduction, 
characterizes all possible triples $\{A_1,A_2,A_3\}$  
which one can obtain in this way, and it is 
a corollary of the following proposition. 

\begin{proposition}\label{propo2}
Let $q_1$, $q_2$, $q_3\in(0,1)$ be given. 
Let $A_1=0$ and $A_2=I$ in $\M^{3\times 3}$. Then there exist 
$S_1,S_2,S_3,A_3\in\M^{3\times 3}$ satisfying 
conditions {\rm \eqref{cond1}} and {\rm \eqref{cond2}} of 
{\rm Lemma \ref{propo1}}, if and only if $A_3$ 
is of the form 
\begin{equation}\label{atre}
A_3=\frac{1}{q_3}
\left[
\left(
1-\prod\limits_{i=1}^3(1-q_i)
\right)
G^{-1}\left(\begin{array}{ccc}
\lambda_1 & 0         & 0         \\
0         & \lambda_2 & 0         \\
0         & 0         & \lambda_3 \\
\end{array}\right)
G-q_2(1-q_3)I
\right]\,,
\end{equation}
where $G$ is an arbitrary matrix in $GL(3)$ and 
the $\lambda_i$'s are defined as follows
$$ 
\lambda_1=0\,,\quad \lambda_2=\frac{1}{1-q_1}\,,\quad 
\lambda_3=\frac{q_2}{q_1+q_2-q_1q_2}\,.
$$
\end{proposition}
\begin{proof}
We rewrite \eqref{cond2} more explicitly:\\[-10pt]
\begin{align}\label{cond2'}
\nonumber & S_2=q_1 A_1+(1-q_1)S_1 =(1-q_1)S_1,\\
          & S_3=q_2 A_2+(1-q_2)S_2 =q_2 I+(1-q_2)(1-q_1)S_1,\\
\nonumber & S_1=q_3 A_3+(1-q_3)S_3 =q_3 A_3+q_2(1-q_3)I+ (1-q_1)(1-q_2)(1-q_3)S_1 \,.
\end{align}
Now let $\lambda_i$ be the eigenvalues of $S_1$, then by \eqref{cond1} and 
\eqref{cond2'} we get \\[-10pt]
\begin{align*}
&\det(A_1-S_1)=0 \Longleftrightarrow \det(-S_1)=0 
\Longleftrightarrow \lambda_1=0\,,\\
&\det(A_2-S_2)=0 \Longleftrightarrow \det(I-(1-q_1)S_1)=0 
\Longleftrightarrow \lambda_2=\frac{1}{1-q_1}\,.
\end{align*}
Moreover one can check that 
$$
\det(A_3-S_3)=0 \Longleftrightarrow \det[(q_1+q_2-q_1q_2)S_1-q_2 I]=0 
\Longleftrightarrow \lambda_3=\frac{q_2}{q_1+q_2-q_1q_2}\,.\\
$$
Note that the $\lambda_i$'s  
are all distinct, since $q_1,q_2,q_3\in(0,1)$. 
Therefore the matrix $S_1$ is diagonalizable. Hence 
for any  
$S_1\in\{G^{-1}{\rm diag}(\lambda_1,\lambda_2,\lambda_3)G,\,\, G\in GL(3)\}$, 
the matrices $A_3$, $S_2$, $S_3$ are uniquely determined by \eqref{cond2'}. 
In particular $A_3$ is of the form \eqref{atre}. Conversely, for any 
$A_3$ of the form $\eqref{atre}$, the matrices $S_1$, $S_2$, $S_3$ are 
uniquely determined and conditions \eqref{cond1} and \eqref{cond2} are satisfied.
\end{proof}

\begin{proof}[{\bf Proof of Theorem \ref{fine}}]
Let $A,M,N\in\M^{3\times 3}$ satisfy the assumptions of Theorem \ref{fine}. 
Without loss of generality we can assume that $M=0$ and $N=I$. 
Notice that $A$ is of the form \eqref{atre}. 
Then, by Proposition \ref{propo2} and Lemma \ref{propo1}, 
the set $K=\{0,I,A\}$ is non-rigid for approximate solutions. 
\end{proof}

\begin{remark}\label{matgn}{\rm 
If we choose $G=I$ in $\eqref{atre}$ the set $K$ reduces itself to 
that given in \cite{GN}, that is 
\begin{align*}
& A_3={\rm diag}{\displaystyle\left(\frac{-q_2(1-q_3)}{q_3}\,,
\frac{q_1+q_3-q_1q_3}{q_3(1-q_1)}\,,
\frac{q_2}{q_1+q_2-q_1q_2}\right)}\,,\\
& S_1={\rm diag}{\displaystyle \left(0,\frac{1}{1-q_1},
\frac{q_2}{q_1+q_2-q_1q_2}\right)}\,,\\
& S_2={\rm diag}{\displaystyle \left(0,1,
\frac{q_2(1-q_2)}{q_1+q_2-q_1q_2}\right)}\,,\\
& S_3={\rm diag}{\displaystyle \left(q_2,1,\frac{q_2}{q_1+q_2-q_1q_2}\right)}\,.
\end{align*}
}
\end{remark}

\medskip
\centerline{\sc Acknowledgements}
The authors wish to thank Enzo Nesi for having proposed the subject 
of this paper, and for many interesting discussions and suggestions.

\end{document}

%% file: div7.pstex_t
\begin{picture}(0,0)%
\includegraphics{div7.pstex}%
\end{picture}%
\setlength{\unitlength}{3039sp}%
\begingroup\makeatletter\ifx\SetFigFont\undefined%
\gdef\SetFigFont#1#2#3#4#5{%
  \reset@font\fontsize{#1}{#2pt}%
  \fontfamily{#3}\fontseries{#4}\fontshape{#5}%
  \selectfont}%
\fi\endgroup%
\begin{picture}(3687,3859)(2101,-4508)
\put(3601,-1711){\makebox(0,0)[lb]{\smash{\SetFigFont{10}{12.0}{\rmdefault}{\mddefault}{\updefault}\special{ps: gsave 0 0 0 setrgbcolor}$Q + w_\delta$\special{ps: grestore}}}}
\put(3601,-4111){\makebox(0,0)[lb]{\smash{\SetFigFont{10}{12.0}{\familydefault}{\mddefault}{\updefault}\special{ps: gsave 0 0 0 setrgbcolor}$O + w_\delta$\special{ps: grestore}}}}
\put(4276,-3736){\makebox(0,0)[lb]{\smash{\SetFigFont{9}{10.8}{\rmdefault}{\mddefault}{\updefault}\special{ps: gsave 0 0 0 setrgbcolor}$r_\delta$\special{ps: grestore}}}}
\put(2851,-1486){\makebox(0,0)[lb]{\smash{\SetFigFont{12}{14.4}{\rmdefault}{\mddefault}{\updefault}\special{ps: gsave 0 0 0 setrgbcolor}$D^2_\delta$\special{ps: grestore}}}}
\put(2851,-3886){\makebox(0,0)[lb]{\smash{\SetFigFont{12}{14.4}{\rmdefault}{\mddefault}{\updefault}\special{ps: gsave 0 0 0 setrgbcolor}$D^1_\delta$\special{ps: grestore}}}}
\put(2101,-2686){\makebox(0,0)[lb]{\smash{\SetFigFont{12}{14.4}{\rmdefault}{\mddefault}{\updefault}\special{ps: gsave 0 0 0 setrgbcolor}$L_\delta$\special{ps: grestore}}}}
\put(3826,-1186){\makebox(0,0)[lb]{\smash{\SetFigFont{10}{12.0}{\rmdefault}{\mddefault}{\updefault}\special{ps: gsave 0 0 0 setrgbcolor}$e_n=v_r$\special{ps: grestore}}}}
\put(5101,-2611){\makebox(0,0)[lb]{\smash{\SetFigFont{10}{12.0}{\rmdefault}{\mddefault}{\updefault}\special{ps: gsave 0 0 0 setrgbcolor}$\nu(s)$\special{ps: grestore}}}}
\end{picture}